\documentclass{agtart_a}
\pdfoutput=1

\usepackage{graphicx} 
\usepackage{pinlabel}

\title{Intrinsic linking and knotting of graphs in arbitrary 3--manifolds}

\author{Erica Flapan}
\givenname{Erica}
\surname{Flapan}
\address{Department of Mathematics\\
Pomona College\\
Claremont, CA 91711\\USA}
\email{eflapan@pomona.edu}
\urladdr{}

\author{Hugh Howards}
\givenname{Hugh}
\surname{Howards}
\address{Department of Mathematics\\
Wake Forest University\\
Winston-Salem, NC 27109\\USA}
\email{howards@wfu.edu}
\urladdr{}

\author{Don Lawrence}
\givenname{Don}
\surname{Lawrence}
\address{Department of Mathematics\\
Occidental College\\
Los Angeles, CA 90041\\USA}
\email{DonL@oxy.edu}
\urladdr{}

\author{Blake Mellor}
\givenname{Blake}
\surname{Mellor}
\address{Department of Mathematics\\
Loyola Marymount University\\
Los Angeles, CA 90045\\USA}
\email{bmellor@lmu.edu}
\urladdr{}

\volumenumber{6}
\issuenumber{}
\publicationyear{2006}
\papernumber{37}
\startpage{1025}
\endpage{1035}

\doi{}
\MR{}
\Zbl{}


\received{25 October 2005}
\revised{3 May 2006}
\accepted{11 May 2006}
\published{9 August 2006}
\publishedonline{9 August 2006}
\proposed{}
\seconded{}
\corresponding{}
\editor{}
\version{}

\arxivreference{math.GT/0508004}

\AtBeginDocument{\let\bar\wbar\let\tilde\wtilde\let\hat\what}

\theoremstyle{plain}
\newtheorem{proposition}{Proposition}
\newtheorem{fact}{Fact}
\newtheorem{theorem}{Theorem}
\newtheorem{lemma}{Lemma}
\theoremstyle{definition}
\newtheorem*{remark}{Remark}

\def\scm{{\wwtilde{M}}} 

\begin{document}

\begin{asciiabstract}
We prove that a graph is intrinsically linked in an arbitrary 3-manifold M if 
and only if it is intrinsically linked in S^3.  Also, assuming the Poincare 
Conjecture, we prove that a graph is intrinsically knotted in M if and only if
it is intrinsically knotted in S^3.
\end{asciiabstract}

\begin{htmlabstract}
We prove that a graph is intrinsically linked in an arbitrary
3&ndash;manifold M if and only if it is intrinsically linked in
S<sup>3</sup>.  Also, assuming the Poincar&eacute; Conjecture, we
prove that a graph is intrinsically knotted in M if and only if it is
intrinsically knotted in S<sup>3</sup>.
\end{htmlabstract}

\begin{abstract}
We prove that a graph is intrinsically linked in an arbitrary 3--manifold $M$ 
if and only if it is intrinsically linked in $S^3$.  Also, assuming the 
Poincar\'e Conjecture, we prove that a graph is intrinsically knotted in $M$ if
and only if it is intrinsically knotted in $S^3$.
\end{abstract}

\maketitle

\section{Introduction}
\label{intro}

The study of intrinsic linking and knotting began in 1983 when Conway and 
Gordon \cite{CG} showed that every embedding of $K_6$ (the complete graph on 
six vertices) in $S^3$ contains a non-trivial link, and every embedding of 
$K_7$ in $S^3$ contains a non-trivial knot.  Since the existence of such a 
non-trivial link or knot depends only on the graph and not on the particular
embedding of the graph in $S^3$, we say that $K_6$ is {\it intrinsically 
linked\/} and $K_7$ is {\it intrinsically knotted\/}.  

\begin{figure}[ht!]
\labellist
\small
\pinlabel $K_6$ at 17 133
\pinlabel $K_{3,3,1}$ at 122 139
\endlabellist
\centering
\includegraphics{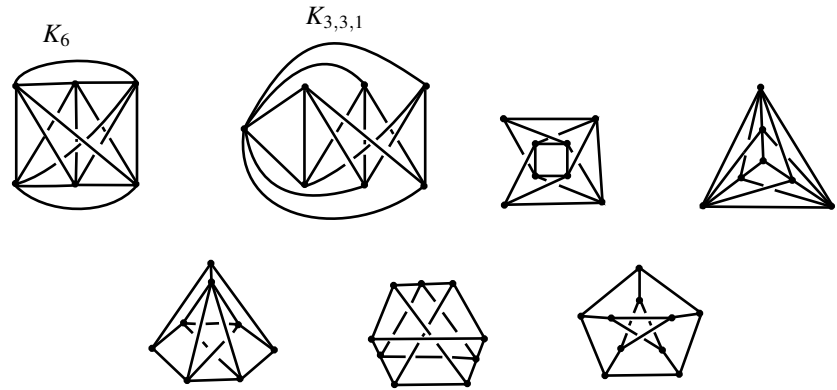}
\caption{The Petersen family of graphs}
\label{petersen}
\end{figure}

At roughly the same time as Conway and Gordon's result, Sachs \cite{S2,S1}
independently proved that $K_6$ and  $K_{3,3,1}$ are intrinsically linked,
and used these two results to prove that any graph with a minor in the
{\it Petersen family\/} (\fullref{petersen}) is intrinsically linked.
Conversely, Sachs conjectured that any graph which is intrinsically linked
contains a minor in the Petersen family.  In 1995, Robertson, Seymour
and Thomas \cite{RST} proved Sachs' conjecture, and thus completely
classified intrinsically linked graphs.

Examples of intrinsically knotted graphs other than $K_7$ are now known, see 
Foisy \cite{F}, Kohara and Suzuki \cite{KS} and Shimabara \cite{Sh}.  
Furthermore, a result of Robertson and Seymour \cite{RS} implies that there are
only finitely many intrinsically knotted graphs that are minor-minimal with 
respect to intrinsic knottedness.  However, as of yet,  intrinsically knotted 
graphs have not been classified.

In this paper we consider the properties of intrinsic linking and knotting in
{\it arbitrary\/} 3--manifolds.  We show that these properties are truly
{\it intrinsic\/} to a graph in the sense that they do not depend on either the
ambient 3--manifold or the particular embedding of the graph in the 
3--manifold.  Our proof in the case of intrinsic knotting assumes the 
Poincar\'e Conjecture.

We will use the following terminology.  By a {\it graph\/} we shall mean a
finite graph, possibly with loops and repeated edges.  Manifolds may have
boundary and do not have to be compact.  All spaces are piecewise linear;
in particular, we assume that the image of an {\it embedding\/} of a
graph in a 3--manifold is a piecewise linear subset of the 3--manifold.
An embedding of a graph $G$ in a 3--manifold $M$ is {\it unknotted\/}
if every circuit in $G$ bounds a disk in $M$; otherwise, the embedding
is {\it knotted\/}.  An embedding of a graph $G$ in a 3--manifold $M$
is {\it unlinked\/} if it is unknotted and every pair of disjoint
circuits in $G$ bounds disjoint disks in $M$; otherwise, the embedding
is {\it linked\/}.  A graph is {\it intrinsically linked\/} in $M$
if every embedding of the graph in $M$ is linked; and a graph is {\it
intrinsically knotted\/} in $M$ if every embedding of the graph in $M$
is knotted.  (So by definition an intrinsically knotted graph must be
intrinsically linked, but not vice-versa.)

The main results of this paper are that a graph is intrinsically linked in an
arbitrary 3--manifold if and only if it is intrinsically linked in $S^3$ 
(\fullref{linked}); and  (assuming the Poincar\'e Conjecture) that a graph
is intrinsically knotted in an arbitrary 3--manifold  if and only if it is 
intrinsically knotted in $S^3$ (\fullref{knotted}).  We use Robertson, 
Seymour, and Thomas' classification of intrinsically linked graphs in $S^3$ for
our proof of \fullref{linked}.  However,  because there is no analogous 
classification of intrinsically knotted graphs in $S^3$, we need to take a 
different approach to prove \fullref{knotted}.  In particular, the proof of
\fullref{knotted} uses \fullref{compact theorem} (every compact 
subset of a simply connected 3--manifold is homeomorphic to a subset of $S^3$),
whose proof in turn relies on the Poincar\'e Conjecture.  Our assumption of the
Poincar\'e Conjecture seems reasonable, because Perelman \cite{P1,P2} has
announced a proof of Thurston's Geometrization Conjecture, which implies the 
Poincar\'e Conjecture \cite{M}.  (See also Morgan and Tian \cite{MT}.)

We would like to thank Waseda University, in Tokyo, for hosting the 
International Workshop on Knots and Links in a Spatial Graph, at which this 
paper was conceived.  We also thank the Japan Society for the Promotion of 
Science for providing funding for the third author with a Grant-in-Aid for 
Scientific Research.

\section{Intrinsically linked graphs}
\label{il}

In this section, we prove that intrinsic linking is independent of the 
3--manifold in which a graph is embedded.  We begin by showing (\fullref{lift})
that any unlinked embedding of a graph $G$ in a 3--manifold lifts to an 
unlinked embedding of $G$ in the universal cover.  In the universal cover, the 
linking number can be used to analyze intrinsic linking (\fullref{edge}), as
in the proofs of Conway and Gordon \cite{CG} and Sachs \cite{S2,S1}.  
After we've shown that $K_6$ and $K_{3,3,1}$ are intrinsically linked in any 
3--manifold (\fullref{k6}), we use the classification of intrinsically 
linked graphs in $S^3$, Robertson, Seymour, and Thomas \cite{RST}, to conclude 
that any graph that is intrinsically linked in $S^3$ is intrinsically linked in
every 3--manifold (\fullref{linked}).

We call a circuit of length 3 in a graph a {\it triangle\/} and a circuit of 
length 4 a {\it square\/}.

\begin{lemma}
\label{lift}
Any unlinked embedding of a graph $G$ in a 3--manifold $M$ lifts to an unlinked
embedding of $G$ in the universal cover $\scm$.
\end{lemma}

\begin{proof}
Let $f\co G \rightarrow M$ be an unlinked embedding.  $\pi_1(G)$ is generated
by the circuits of $G$ (attached to a basepoint).  Since $f(G)$ is
unknotted, every cycle in $f(G)$ bounds a disk in $M$.  So $f_*(\pi_1(G))$
is trivial in $\pi_1(M)$.

Thus, an unlinked embedding of $G$ into $M$ lifts to an embedding of $G$
in the universal cover $\scm$.  Since the embedding into $M$ is unlinked,
cycles of $G$ bound disks in $M$ and pairs of disjoint cycles of $G$
bound disjoint disks in $M$.  All of these disks in $M$ lift to disks
in $\scm$, so the embedding of the graph in $\scm$ is also unlinked.
\end{proof}

Recall that if $M$ is a 3--manifold with $H_1(M)=0$, then disjoint oriented 
loops $J$ and $K$ in $M$ have a well-defined linking number $\mathrm{lk}(J,K)$,
which is the algebraic intersection number of $J$ with any oriented surface 
bounded by $K$.  Also, the linking number is symmetric:  
$\mathrm{lk}(J,K)=\mathrm{lk}(K,J)$.  

It will be convenient to have a notation for the linking number modulo 2:  
Define $\omega(J,K)=\mathrm{lk}(J,K)\bmod{2}$.  Notice that $\omega(J,K)$ is 
defined for a pair of {\it unoriented\/} loops.  Since linking number is 
symmetric, so is $\omega(J,K)$.  If $J_1$, \ldots, $J_n$ are loops in an 
embedded graph such that in the list $J_1,\ldots,J_n$ every edge appears an 
even number of times, and if $K$ is another loop, disjoint from the $J_i$, then
$\sum \omega(J_i,K)=0\bmod{2}$.

If $G$ is a graph embedded in a simply connected 3--manifold, let 
$$\omega(G)=\sum \omega(J,K)\bmod{2},$$
where the sum is taken over all {\it unordered\/} pairs $(J,K)$ of
disjoint circuits in $G$.  Notice that if $\omega(G)\neq 0$, then the
embedding is linked (but the converse is not true).

\begin{lemma}
\label{edge}
Let $\scm$ be a simply connected 3--manifold, and let $H$ be an embedding of 
$K_6$ or $K_{3,3,1}$ in $\scm$.  Let $e$ be an edge of $H$, and let $e'$ be an
arc in $\scm$ with the same endpoints as $e$, but otherwise disjoint from $H$.
Let $H'$ be the graph $(H-e)\cup e'$.  Then $\omega(H')=\omega(H)$.
\end{lemma}

\begin{proof}
Let $D=e\cup e'$.

First consider the case that $H$ is an embedding of $K_6$.  We will count how
many terms in the sum defining $\omega(H)$ change when $e$ is replaced by $e'$.
Let $K_1$, $K_2$, $K_3$ and $K_4$ be the four triangles in $H$ disjoint from 
$e$ (hence also disjoint from $e'$ in $H'$), and for each $i$ let $J_i$ be the
triangle complementary to $K_i$.  The $J_i$ all contain $e$.  For each $i$, let
$J'_i=(J_i-e)\cup e'$, and notice that
\begin{equation} \label{eq1}
\omega(J'_i,K_i)=\omega(J_i,K_i)+\omega(D,K_i)\bmod{2}.
\end{equation}

Because each edge appears twice in the list $K_1, K_2, K_3, K_4$, we have 
$\omega(K_1,D)+\omega(K_2,D)+\omega(K_3,D)+\omega(K_4,D)=0\bmod{2}$.  Thus, 
$\omega(K_i,D)$ is nonzero for an even number of $i$.  It follows from 
\fullref{eq1} that there are an even number of $i$ such that 
$\omega(J'_i,K_i)\neq\omega(J_i,K_i)$.  Thus, 
$\sum_{i=1}^4\omega(J'_i,K_i)=\sum_{i=1}^4\omega(J_i,K_i)\bmod{2}$, and
\begin{eqnarray*}
 \omega(H') & = & \sum_{\stackrel{J,K\subseteq H'}{\ni e'\notin J,K}} \omega(J,K)
                  +\sum_{i=1}^4\omega(J'_i,K_i)\bmod{2} \\
            & = & \sum_{\stackrel{J,K\subseteq H}{\ni e\notin J,K}} \omega(J,K)
                  +\sum_{i=1}^4\omega(J_i,K_i)\bmod{2} \\
            & = & \omega(H)
\end{eqnarray*}	

Next consider the case that $H$ is an embedding of $K_{3,3,1}$.  Let $x$ be the
vertex of valence six in $H$ (and in $H'$).

{\bf Case 1}\qua $e$ contains $x$.  Then $e$ is not in any square in $H$ that 
has a complementary disjoint triangle.  Let $K_1$, $K_2$ and $K_3$ be the three
squares in $H$ disjoint from $e$, and let $J_1$, $J_2$ and $J_3$ be the
corresponding complementary triangles, all of which contain $e$.  As in the 
$K_6$ case, let $J'_i=(J_i-e)\cup e'$ for each $i$; again we have 
\fullref{eq1}.  Every edge in the list $K_1, K_2, K_3$ appears exactly twice, 
so $\omega(K_1,D)+\omega(K_2,D)+\omega(K_3,D)=0\bmod{2}$.  Thus, 
$\omega(K_i,D)$ is nonzero for an even number of $i$; and for an even number of
$i$, $\omega(J'_i,K_i)\neq\omega(J_i,K_i)$.  The other pairs of circuits 
contributing to $\omega(H)$ do not involve $e$.  As in the $K_6$ case, it 
follows that $\omega(H')=\omega(H)$.

{\bf Case 2}\qua $e$ doesn't contain $x$.  Let $J_0$ be the triangle containing
$e$, and let $K_0$ be the complementary square.  Let $J_1$ through $J_4$ be the
four squares that contain $e$, but not $x$ (so that they have complementary 
triangles); and let $K_1$ through $K_4$ be the complementary triangles.  With 
$J'_i$ defined as in the other cases, we again have \fullref{eq1}.  Every edge 
appears an even number of times in the list $K_0, K_1, K_2, K_3, K_4$, so  
$\sum_{i=0}^4\omega(K_i,D)=0\bmod{2}$, and $\omega(K_i,D)\neq0$ for an even
number of $i$.  As in the other cases, it follows that for an even number of
$i$, $\omega(J'_i,K_i)\neq\omega(J_i,K_i)$; and an even number of the terms in
the sum defining $\omega(H)$ change when $e$ is replaced by $e'$; and 
$\omega(H')=\omega(H)$.
\end{proof}

\begin{proposition}
\label{k6}
$K_6$ and $K_{3,3,1}$ are intrinsically linked in any 3--manifold $M$.
\end{proposition}

\begin{proof}
Let $G$ be either $K_6$ or $K_{3,3,1}$, and let $f\co G\rightarrow M$ be an 
embedding. Suppose for the sake of contradiction that $f(G)$ is unlinked.  Let
$\scm$ be the universal cover of $M$.  By \fullref{lift}, $f$ lifts to an 
unlinked embedding $\wtilde{f}\co G\rightarrow\scm$.

Let $\wtilde{G}=\wtilde{f}(G)\subseteq\scm$, and let $\wtilde{H}$ be a
copy of $G$ embedded in a ball in $\scm$.  Isotope $\wtilde{G}$ so that 
$\wtilde{H}$ and $\wtilde{G}$ have the same vertices, but do not 
otherwise intersect.  Then $\wtilde{G}$ can be transformed into 
$\wtilde{H}$ by changing one edge at a time -- replace an edge of 
$\wtilde{G}$ by the corresponding edge of $\wtilde{H}$, once for every 
edge.  By repeated applications of \fullref{edge}, 
$\omega(\wtilde{G})=\omega(\wtilde{H})$.  Since $\wtilde{H}$ is inside
a ball in $\scm$, Conway and Gordon's proof \cite{CG}, and Sachs' proof 
\cite{S2,S1}, that $K_6$ and $K_{3,3,1}$ are intrinsically linked in 
$S^3$, show that $\omega(\wtilde{H})=1$. 

Thus, $\omega(\wtilde{G})=1$, and there must be disjoint circuits $J$ and 
$K$ in $\wtilde{G}$ that do not bound disjoint disks in $\scm$, 
contradicting that $\tilde{f}$ is an {\it unlinked} embedding.  Thus, $f(G)$ is
linked in $M$.
\end{proof}

\begin{figure}[ht!]
\labellist
\pinlabel $a$ [t] at 53 19
\pinlabel $b$ [br] at 10 96
\pinlabel $c$ [bl] at 94 96
\pinlabel $a$ [t] at 269 17
\pinlabel $b$ [br] at 226 96
\pinlabel $c$ [bl] at 315 97
\pinlabel $v$ [tr] at 267 70
\endlabellist
\centering
\includegraphics{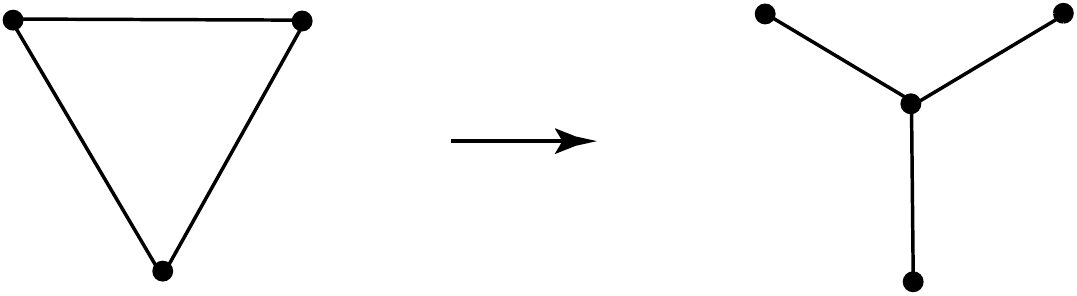}
\caption{A $\Delta-Y$ Move}
\label{triangley}
\end{figure}
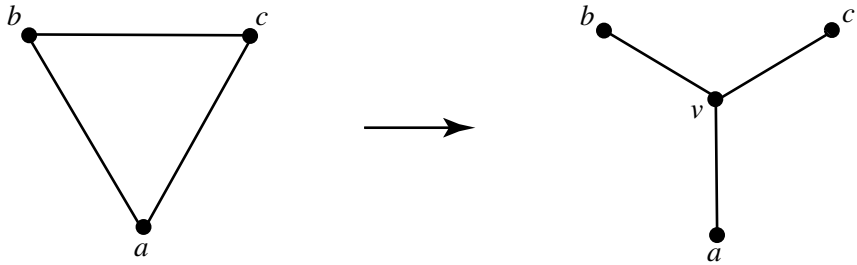

Let $G$ be a graph which contains a triangle $\Delta$.  Remove the three edges 
of $\Delta$ from $G$.  Add three new edges, connecting the three vertices of 
$\Delta$ to a new vertex.  The resulting graph, $G'$, is said to have been 
obtained from $G$ by a ``$\Delta-Y$ move'' (\fullref{triangley}).  The seven
graphs that can be obtained from $K_6$ and $K_{3,3,1}$ by $\Delta-Y$ moves are 
the {\it Petersen family\/} of graphs (\fullref{petersen}).  

If a graph $G'$ can be obtained from a graph $G$ by repeatedly deleting edges 
and isolated vertices of $G$, and/or contracting edges of $G$, then $G'$ is a
{\it minor\/} of $G$.

The following facts were first proved, in the $S^3$ case, by Motwani, 
Raghunathan and Saran \cite{MRS}.  Here we generalize the proofs to any 
3--manifold $M$.

\begin{fact} \label{fact1}
If a graph $G$ is intrinsically linked in $M$, and $G'$ is obtained from $G$ by
a $\Delta-Y$ move, then $G'$ is intrinsically linked in $M$.
\end{fact}

\begin{proof}
Suppose to the contrary that $G'$ has an unlinked embedding $f\co G'\rightarrow
M$.  Let $a$, $b$, $c$ and $v$ be the embedded vertices of the $Y$ 
illustrated in \fullref{triangley}.   Let $B$ denote a regular neighborhood 
of the embedded $Y$ such that $a$, $b$ and $c$ are on the boundary of $B$, $v$ 
is in the interior of $B$, and $B$ is otherwise disjoint from $f(G')$.  Now add
edges $ab$, $bc$ and $ac$ in the boundary of $B$ so that the resulting 
embedding of the $K_4$ with vertices $a$, $b$, $c$, and $v$ is panelled in $B$
(ie, every cycle bounds a disk in the complement of the graph).  We now remove 
vertex $v$ (and its incident edges) to get an embedding $h$ of $G$ such that if
$e$ is any edge of $G\cap G'$ then $h(e)=f(e)$ and the triangle $abc$ is in 
$\partial B$.  

Observe that if $K$ is any circuit in $h(G)$ other than the triangle $abc$, 
then $K$ is isotopic to a circuit in $G'$.  The triangle $abc$ bounds a disk in
$B$, and since $f(G')$ is unknotted, every circuit in $f(G')$ bounds a disk in
$M$.  Thus $h(G)$ is unknotted.  Also if $J$ and $K$ are disjoint circuits in 
$h(G)$ neither of which is $abc$, then $J\cup K$ is isotopic to a pair of 
disjoint circuits $J'\cup K'$ in $f(G')$.  Since $f(G')$ is unlinked, $J'$ and 
$K'$ bound disjoint disks in $M$.  Hence $J$ and $K$ also bound disjoint disks 
in $M$. Finally if $K$ is a circuit in $h(G)$ which is disjoint from $abc$, 
then $K$ is contained in $f(G')$.  Since $f(G')$ is unknotted, $K$ bounds a 
disk $D$ in $M$.  Furthermore, since $B$ is a ball,  we can isotope $D$ to a 
disk which is disjoint from $B$.  Now $abc$ and $K$ bound disjoint disks in 
$M$.  So $h(G)$ is unlinked, contradicting the hypothesis that $G$ is 
intrinsically linked in $M$.  We conclude that $G'$ is also intrinsically 
linked in $M$.
\end{proof}

\begin{fact} \label{fact2}
If a graph $G$ has an unlinked embedding in $M$, then so does every minor of 
$G$.
\end{fact}
\begin{proof}
The proof is identical to the proof for $S^3$.
\end{proof}

\begin{theorem}
\label{linked}
Let $G$ be a graph, and let $M$ be a 3--manifold.  The following are 
equivalent:
\begin{itemize}
  \item[(1)] $G$ is intrinsically linked in $M$,
  \item[(2)] $G$ is intrinsically linked in $S^3$,
  \item[(3)] $G$ has a minor in the Petersen family of graphs.
\end{itemize}
\end{theorem}

\begin{proof}
Robertson, Seymour and Thomas \cite{RST} proved that (2) and (3) are 
equivalent.  We see as follows that (1) implies (2):  Suppose there is an 
unlinked embedding of $G$ in $S^3$.  Then the embedded graph and its system of 
disks in $S^3$ are contained in a ball, 
which embeds in $M$.

We will complete the proof by checking that (3) implies (1).   $K_6$ and 
$K_{3,3,1}$ are intrinsically linked in $M$ by \fullref{k6}.  Thus, by 
\fullref{fact1}, all the graphs in the Petersen family are intrinsically linked
in $M$.  Therefore, if $G$ has a minor in the Petersen family, then it is 
intrinsically linked in $M$, by \fullref{fact2}.
\end{proof}

\section{Compact subsets of a simply connected space}
\label{compact}

In this section, we assume the Poincar\'e Conjecture, and present some known 
results about 3--manifolds, which will be used in \fullref{ik} to prove 
that intrinsic knotting is independent of the 3--manifold (\fullref{knotted}).
  
\begin{fact}
\label{balls}
Assume that the Poincar\'e Conjecture is true.  Let $\scm$ be a simply 
connected 3--manifold, and suppose that $B\subseteq\scm$ is a compact 
3--manifold whose boundary is a disjoint union of spheres.  Then $B$ is a ball 
with holes (possibly zero holes).
\end{fact}

\begin{proof}
By the Seifert--Van Kampen theorem, $B$ itself is simply connected.  Cap off 
each boundary component of $B$ with a ball, and the result is a closed simply
connected 3--manifold.  By the Poincar\'e Conjecture, this must be the 
3--sphere.
\end{proof}

\begin{fact}
\label{compress}
Let $\scm$ be a simply connected 3--manifold, and suppose that $N\subseteq\scm$
is a compact 3--manifold whose boundary is nonempty and not a union of spheres.
Then there is a compression disk $D$  in $\scm$ for a component of $\partial N$
such that $D\cap\partial N=\partial D$.
\end{fact}

\begin{proof}
Since $\partial N$ is nonempty, and not a union of spheres, there is a boundary
component $F$ with positive genus.  Because $\scm$ is simply connected, $F$ is 
not incompressible in $\scm$.  Thus, $F$ has a compression disk.

Among all compression disks for boundary components of $N$ (intersecting 
$\partial N$ transversely), let $D$ be one such that $D\cap\partial N$ consists
of the fewest circles.  Suppose, for the sake of contradiction, that there is a
circle of intersection in the interior of $D$.  Let $c$ be a circle of 
intersection which is innermost in $D$, bounding a disk $D'$ in $D$.  Either 
$c$ is nontrivial in $\pi_1(\partial N)$, in which case $D'$ is itself a 
compression disk; or $c$ is trivial, bounding a disk on $\partial N$, which can
be used to remove the circle $c$ of intersection from $D\cap\partial N$.  In 
either case, there is a compression disk for $\partial N$ which has fewer 
intersections with $\partial N$ than $D$ has, contradicting minimality.  Thus, 
$D\cap\partial N=\partial D$.
\end{proof}
  
We are now ready to prove the main result of this section.  Because its proof 
uses \fullref{balls}, it relies on the Poincar\'e Conjecture.
\medskip

\begin{proposition}
\label{compact theorem}
Assume that the Poincar\'e Conjecture is true.  Then every compact subset $K$
of a simply connected 3--manifold $\scm$ is homeomorphic to a subset of $S^3$.
\end{proposition}    

\begin{proof}
We may assume without loss of generality that $K$ is connected.  Let 
$N\subseteq\scm$ be a closed regular neighborhood of $K$ in $\scm$.  Then $N$
is a compact connected 3--manifold with boundary.  It suffices to show that $N$
embeds in $S^3$.

Let $g(S)$ denote the genus of a connected closed orientable surface $S$.  
Define the complexity $c(S)$ of a closed orientable surface $S$ to be the sum 
of the squares of the genera of the components $S_i$ of $S$, so 
$c(S) = \sum_{S_i}{g(S_i)^2}$.  Our proof will proceed by induction on 
$c(\partial N)$.  We make two observations about the complexity function.
\begin{enumerate}
  \item $c(S) = 0$ if and only if $S$ is a union of spheres.
  \item If $S'$ is obtained from $S$ by surgery along a non-trivial simple 
        closed curve $\gamma$, then $c(S') < c(S)$.  
\end{enumerate}
	
We prove Observation (2) as follows.  It is enough to consider the 
component $S_0$ of $S$ containing $\gamma$.  If $\gamma$ separates $S_0$, then 
$S_0 = S_1 \# S_2$, where $S_1$ and $S_2$ are not spheres, and $S'$ is the 
result of replacing $S_0$ by $S_1 \cup S_2$ in $S$.  In this case, 
$c(S_0) = g(S_0)^2 = (g(S_1) + g(S_2))^2 = c(S_1) + c(S_2) + 2g(S_1)g(S_2) > 
c(S_1) + c(S_2)$, since $g(S_1)$ and $g(S_2)$ are nonzero.  On the other hand, 
if $\gamma$ does not separate $S_0$, then surgery along $\gamma$ reduces the 
genus of the surface.  Then the square of the genus is also smaller, and hence 
again $c(S') < c(S)$.

If $c(\partial N) = 0$, then by \fullref{balls} $N$ is a ball with holes, and 
so embeds in $S^3$, establishing our base case.  If $c(\partial N) > 0$, then 
by \fullref{compress} there is a compression disk $D$ for $\partial N$ such 
that $D\cap\partial N=\partial D$.  There are three cases to consider.

{\bf Case 1}\qua $D\cap N=\partial D$.  Let $N' = N \cup {\rm nbd}(D)$.  Since 
$\partial N'$ is the result of surgery on $\partial N$ along a non-trivial 
simple closed curve, $c(\partial N') < c(\partial N)$, so by induction $N'$ 
embeds in $S^3$.  Hence $N$ embeds in $S^3$.

{\bf Case 2}\qua $D\cap N= D$, and $D$ separates $N$.  Then cutting $N$ along 
$D$ (ie removing $D \times (-1, 1)$) yields two connected manifolds $N_1$ and
$N_2$, with $c(\partial N_1) < c(\partial N)$ and $c(\partial N_2) < 
c(\partial N)$.  So $N_1$ and $N_2$ each embed in $S^3$.  Consider two copies 
of $S^3$, one containing $N_1$ and the other containing $N_2$.

Let $C_1$ be the component of $S^3 - N_1$ whose boundary contains 
$D \times \{1\}$, and $C_2$ be the component of $S^3 - N_2$ whose boundary 
contains $D \times \{-1\}$.
Remove small balls $B_1$ and $B_2$ from $C_1$ and $C_2$, respectively.
Then glue together the balls $\mathrm{cl}(S^3-B_1)$ and 
$\mathrm{cl}(S^3-B_2)$ along their boundaries.  The result is a 3--sphere 
containing both $N_1$ and $N_2$, in which $D \times \{1\}$ and $D \times 
\{-1\}$ lie in the boundary of the same component of $S^3 - (N_1 \cup N_2)$. 
So we can embed the arc $\{0\} \times (-1, 1)$ (the core of $D \times (-1, 1)$)
in $S^3 - (N_1 \cup N_2)$, which means we can extend the embedding of 
$N_1 \cup N_2$ to an embedding of $N$.

{\bf Case 3}\qua $D\cap N= D$, but $D$ does not separate $N$.  Then cutting $N$
along $D$ yields a new connected manifold $N'$ with $c(\partial N') < 
c(\partial N)$, so $N'$ embeds in $S^3$.  As in the last case, we also need to 
embed the core $\gamma$ of $D$.  Suppose for the sake of contradiction that 
$\gamma$ has endpoints on two different boundary components $F_1$ and $F_2$ of
$N'$.  Let $\beta$ be a properly embedded arc in $N'$ connecting $F_1$ and 
$F_2$.  Then $\gamma\cup\beta$ is a loop in $\scm$ that intersects the closed 
surface $F_1$ in exactly one point.  But because $H_1(\scm)=0$, the algebraic 
intersection number of $\gamma\cup\beta$ with $F_1$ is zero.  This is 
impossible since $\gamma\cup\beta$ meets $F_1$ in a single point.  Thus, both 
endpoints of $\gamma$ lie on the same boundary component of $N'$, and so 
$\gamma$ can be embedded in $S^3 - N'$.  So the embedding of $N'$ can be 
extended to an embedding of $N$ in $S^3$.
\end{proof}
                 
\section{Intrinsically knotted graphs}
\label{ik}

In this section, we use \fullref{compact theorem} to prove that the 
property of a graph being intrinsically knotted is independent of the 
3--manifold it is embedded in.  Notice that since \fullref{compact theorem} 
relies on the Poincar\'e Conjecture, so does the intrinsic knotting result.

\begin{theorem}
\label{knotted}
Assume that the Poincar\'e Conjecture is true.  Let $M$ be a 3--manifold.  A 
graph is intrinsically knotted in $M$ if and only if it is intrinsically 
knotted in $S^3$.
\end{theorem}
\begin{proof}
Suppose that a graph $G$ is not intrinsically knotted in $S^3$.  Then it embeds
in $S^3$ in such a way that every circuit bounds a disk embedded in $S^3$.  The
union of the embedding of $G$ with these disks is compact, hence is contained 
in a ball $B$ in $S^3$.  Any embedding of $B$ in $M$ yields an unknotted 
embedding of $G$ in $M$.

Conversely, suppose there is an unknotted embedding $f\co G\rightarrow M$.
Let $\scm$ be the universal cover of $M$. By using the same argument
as in the proof of Lemma 1, we can lift $f$ to an unknotted embedding
$\tilde{f}:G\rightarrow\scm$.  Let $K$ be the union of $\tilde{f}(G)$
with the disks bounded by its circuits.   Then $K$ is compact, so by
Proposition \ref{compact theorem}, there is an embedding $g:K\rightarrow
S^3$.  Now $g\circ\tilde{f}(G)$ is an embedding of $G$ in $S^3$, in
which every circuit bounds a disk.  Hence $g\circ\tilde{f}(G)$ is an
unknotted embedding of $G$ in $S^3$.
\end{proof}

\begin{remark}
The proof of \fullref{knotted} can also be used, almost verbatim, to show 
that intrinsic {\it linking\/} is independent of the 3--manifold.
Of course, this argument relies on the Poincar\'e Conjecture; so the
proof given in \fullref{il} is more elementary.
\end{remark}

\bibliographystyle{gtart}
\bibliography{link}

\end{document}